\UseRawInputEncoding
\documentclass[letterpaper, 10 pt, conference]{ieeeconf}
\IEEEoverridecommandlockouts                       
\usepackage{amsmath}    
\usepackage{amssymb}  
\usepackage{graphicx} 
\usepackage{dblfloatfix}    
\usepackage{multirow}
\usepackage{eqnarray,balance}
\usepackage[compress,sort]{cite}
\usepackage{algorithm}
\usepackage[noend]{algpseudocode}

\title{\LARGE \bf
Surrogate Neural Networks for Efficient Simulation-based Trajectory Planning Optimization}

\author{Evelyn Ruff, Rebecca Russell, Matthew Stoeckle, Piero Miotto, and Jonathan P. How
\thanks{Work supported by the Charles Stark Draper Laboratory, Inc}
\thanks{E.\ Ruff is with the Department of Aeronautics and Astronautics, Massachusetts Institute of Technology, and is a Draper Scholar, Cambridge, MA 02139, USA
        {\tt everuff@mit.edu}}%
\thanks{R.\ Russell, M.\ Stoeckle, \& P.\ Miotto are with the Charles Stark Draper Laboratory, Inc, Cambridge, MA 02139, USA
        {\tt \{pmiotto,mstoeckle,rrussell\}@draper.com}}%
\thanks{J.\ How is the R.\ C.\ Maclaurin Professor of Aeronautics and Astronautics, Massachusetts Institute of Technology, {\tt jhow@mit.edu}}%
}

\begin{document}

\maketitle
\thispagestyle{empty}
\pagestyle{empty}
\setlength{\abovedisplayskip}{3pt}
\setlength{\belowdisplayskip}{3pt}

\begin{abstract}

This paper presents a novel methodology that uses surrogate models in the form of neural networks to reduce the computation time of simulation-based optimization of a reference trajectory. Simulation-based optimization is necessary when there is no analytical form of the system accessible, only input-output data that can be used to create a surrogate model of the simulation. Like many high-fidelity simulations, this trajectory planning simulation is very nonlinear and computationally expensive, making it challenging to optimize iteratively. Through gradient descent optimization, our approach finds the optimal reference trajectory for landing a hypersonic vehicle. 
In contrast to the large datasets used to create the surrogate models in prior literature, our methodology is specifically designed to minimize the number of simulation executions required by the gradient descent optimizer. We demonstrated this methodology to be more efficient than the standard practice of hand-tuning the inputs through trial-and-error or randomly sampling the input parameter space. Due to the intelligently selected input values to the simulation, our approach yields better simulation outcomes that are achieved more rapidly and to a higher degree of accuracy. 
Optimizing the hypersonic vehicle's reference trajectory is very challenging due to the simulation's extreme nonlinearity, but even so, this novel approach found a 74\% better-performing reference trajectory compared to nominal, and the numerical results clearly show a substantial reduction in computation time for designing future trajectories. 

\end{abstract}

\section{INTRODUCTION}

High-fidelity simulations are used to analyze the dynamics of complex systems in many engineering and scientific disciplines. For most of these applications, there is some desired system outcome, such as the curvature of an airfoil or the path planned for a robot, but, lacking a closed-form model, the simulation must be used in the design process. Here arises the need for simulation-based optimization (SO), as opposed to traditional optimization with analytical models, generally with the goal of optimizing the simulation input values to achieve the desired output. 

There are many methods of simulation-based optimization \cite{ammeri,amaran} applied across disciplines. A key issue in these approaches is that the execution of such simulations often requires a large amount of computational power and/or processing time. When the simulation itself is too computationally expensive, surrogate models are needed to minimize the number of queries and, therefore, the total computation time. The surrogate model is designed to imitate the input-output behavior of the simulation. Input-output data is required from the simulation to create this surrogate model, so efficiency in creating that model is a key concern \cite{cozad}. The methodology presented herein uses neural networks as surrogate models and provides a novel sampling strategy to reduce the number of simulation runs -- a crucial concern when each run of the simulation takes several minutes, and optimization might require hundreds, if not thousands \cite{du}, of simulation runs.

The use case simulation for this methodology is an Approach and Landing (A/L) simulation for hypersonic vehicles (HV) combined with the reference trajectory calculation \cite{barton}. The HV simulation includes high-fidelity models for flex, slosh, engine thrust, aerodynamics, and full flight software in addition to the trajectory planning algorithm which propagates the vehicle down many possible trajectories. The simulation takes 2-3 minutes to run and the relationship between inputs and outputs is highly nonlinear due to the series of models listed. The system inputs are the 13 trajectory design parameters and the three outputs are the most important performance metrics of the vehicle during A/L.

Hence, finding the set of trajectory design parameters that provides the optimal landing trajectory for the HV is typically done  by manually tuning the various input parameters, running the full simulation, and evaluating the performance outputs. A Monte Carlo approach could also be implemented by randomly selecting hundred or thousands of input values until a sufficiently good reference trajectory is found. Finding a sufficiently good A/L reference trajectory using either of these approaches is a very time-consuming iterative process that requires extensive domain knowledge to either understand the complicated input-output mapping or to accurately evaluate the performance outcomes.

The algorithm outlined in this paper offers a novel simulation-based optimization approach that eliminates the need for hand-tuning input parameters and significantly reduces the overall computation time. The main contributions of this methodology are:
\begin{itemize}
\item Novel algorithm using neural network gradients to intelligently select simulation input values that improves desired simulation outcome by 74\%.
\item Computationally efficient methodology for optimizing general black-box simulations with desired outputs  shown to be six times faster than Monte Carlo approach;
\item New automated optimal trajectory planning tool for hypersonic vehicles that takes hours as opposed to days.
\end{itemize}

\section{RELATED WORK}

This paper aims to improve upon two categories of research: simulation-based optimization via surrogate models and optimal trajectory planning.

\subsection{Simulation-based Optimization via Surrogate Models}

There are many simulation-based optimization algorithms, and their suitability relies heavily on the specific application. Whether the system is continuous or discrete, cheap or expensive to evaluate, and deterministic or stochastic all must be considered \cite{amaran}. This makes various SO methods very challenging to compare. In literature, these methods are often divided into four categories: statistical selection, surrogate models (metamodels), stochastic gradient estimation, and global search \cite{ammeri}. Surrogate models reduce the computational burden of optimizing the simulation outputs by creating a simpler, less-expensive model of the simulation to use instead of the simulation itself. 

For the purpose of this paper, only SO via surrogate models will be analyzed because the chosen application is continuous and very expensive to run, making surrogate modeling the ideal type of algorithm to use \cite{ammeri}. There have been many kinds of surrogate models applied across various fields. Most notable include algebraic models \cite{cozad}, kriging \cite{simpson}, polynomial response surfaces \cite{li}, and radial basis functions \cite{jakobsson}. However, with the proliferation of machine learning in recent decades, surrogate models have more recently been created using support-vector machines \cite{wan} and neural networks \cite{hurrion}. 

Once the surrogate model has been fitted to the data, it is integrated with an optimizer to find the optimal set of inputs that achieve the desired output to the simulation. The applications where neural networks have been applied as surrogate models are very diverse, such as fluids systems \cite{majumder,sreekanth,mirghani} and airfoil design \cite{sun,tao}. Instead of taking advantage of the gradient information available from the neural network, previous methods generally have relied on gradient-free optimization, such as genetic algorithms \cite{sreekanth,mirghani,sun}, particle swarm optimization \cite{tao}, and the grey wolf optimizer \cite{majumder}.

Ref.~\cite{white} uses the gradients of a neural network in combination with the IPOPT software library to optimize a topology. However, the simulation was simple enough that only a single-layer feedforward neural network was necessary, which meant the derivatives had to be calculated analytically and fed into the optimizer. Ref.~\cite{du} uses a similar method of gradient-based optimization using multi-layer neural networks combined with the SNOPT software library to optimize airfoils. Both methods rely on creating a highly accurate surrogate model that only needs to be trained once to be used by the optimizer. Ref.~\cite{white} does this through a complicated training procedure using the Sobolov error and \cite{du} does this by collecting tens of thousands of data points. In the method proposed in this paper, the algorithm collects more data from the actual simulation iteratively in areas of interest suggested by the surrogate model to reduce the total number of queries and save on computation cost.

\subsection{Optimal Trajectory Planning}

Reference trajectory optimization is a very challenging process that requires formulating a nonlinear optimal control problem with multiple constraints and solving, directly or indirectly \cite{betts}. In general, reference trajectories are calculated offline and stored in the HV's computer, although some algorithms work to reduce computation time in order to be calculated during flight \cite{liu}. For comparison, only offline algorithms will be considered in this review. 

Ref.~\cite{dong} proposes a methodology for HV trajectory optimization by combining particle swarm optimization and non-intrusive polynomial chaos. This method works to minimize the flight time of the HV over the course of the trajectory while being robust to physical uncertainties. While a time-based objective function works for the re-entry phase of hypersonic flight, it does not allow for optimizing specific properties at phase completion. Furthermore, this is a direct numerical method that requires access to the dynamical functions, which is not the case in the application of the algorithm proposed in this paper. The same is true for the generalized polynomial chaos algorithm proposed by \cite{okamoto} to optimize landing trajectories of airplanes and the mapped Chebyshev pseudospectral method presented by \cite{wang}. 

Currently, there are no methods of optimal trajectory planning that use surrogate models for simulation-based optimization.

\section{PRELIMINARIES}

\subsection{Neural Networks as Surrogate Models} 

Neural networks (NNs) are very good at modeling complex relationships between inputs and outputs. In fact, multilayer feedforward networks can approximate any measurable function to any degree of accuracy, provided enough training data \cite{hornik}. Multilayer NNs can find intricate structures in high-dimensional data and learn hierarchical feature representations with multiple levels of abstraction \cite{yann}.
Furthermore, neural networks can simultaneously estimate the derivatives of an approximated function, even if the function does not have classically differentiable functions \cite{hornik}. This makes them a perfect choice for surrogate models. Most commonly, the derivatives of the neural network, referred to as gradients, are used to backpropagate error through the neural network to train them on an initial dataset. These gradients can also be used for sensitivity analysis \cite{white} or for gradient-based optimization, as in this proposed method.

\subsection{Objective Function with Constraints}

For all simulations, there must be $m$-dimensional inputs $\left( x \in \mathbb{R}^{\mathrm{m}} \right)$ into the simulation ${f}( x )$ and $n$-dimensional outputs $\left( y \in \mathbb{R}^{\mathrm{n}} \right)$. Furthermore, there are $m$-dimensional minimum and maximum bounds $\left( x_{\min }, x_{\max } \in \mathbb{R}^{\mathrm{m}}\right)$ for each input as well as a $n$-dimensional desired set of outputs $\left( y_{\text {target }} \in \mathbb{R}^{\mathrm{n}}\right)$.
Using an initial training dataset, a neural network is trained to be a surrogate model $\hat{f}(x)$ for the real simulation. This surrogate model is able to predict $n$-dimensional outputs based on real $m$-dimensional inputs \begin{equation} \hat{y}=\hat{f}(x) \end{equation} where $\hat{y}$ is the predicted outputs from the surrogate model. This surrogate model is then used to minimize the loss between the desired outputs and the predicted outputs from the surrogate model. This loss is evaluated by the objective  
\begin{equation} g\left( \hat{y}, y_{\text {target }}\right). \end{equation}
Note that the parameters being optimized are not the predicted output from the surrogate model, but the inputs to the surrogate model, making the objective a function of $ x $:
\begin{equation} g\left(\hat{f}(x), y_{\text {target }}\right) \end{equation}
Further note that this is a constrained optimization problem because the inputs must be bounded within the initial training data range or else the surrogate model will not be able to accurately predict the input's respective output. The inputs are bounded in the physically feasible range by the inequality conditions:
\vspace{-0.1in}
\begin{eqnarray}
c_1(x): x \geq x_{\min } \\
c_2(x): x \leq x_{\max } 
\vspace{-0.1in}
\end{eqnarray}
Because the objective function is a function of the inputs, its gradients can be used to optimize the inputs to minimize the objective loss. The gradient of the objective function with respect to the inputs is \begin{equation}\frac{\partial}{\partial x}\left[ g(\hat{f}(x)) \right] = g^{\prime}(\hat{f}(x))  \hat{f}^{\prime}(x). \end{equation} 
This derivation shows the need for the surrogate model and its gradients for optimizing the inputs to minimize the objective function.

\section{APPROACH}

This section details the general methodology for using surrogate neural networks to optimize a black-box simulation, as shown in Algorithm 1. Subsections B-H detail how to apply the methodology to a specific trajectory planning simulation.

\vspace{-0.1in}
\begin{algorithm}
\caption{}
\begin{algorithmic}[1]
    \State Query an initial dataset $\mathcal{D}=\left\{x \in \mathbb{R}^{\mathrm{m}}\right\}$ from a quasi-random sample. [\ref{ssec:initdata}]
    \State Train a neural network surrogate model $ \hat{f}(x) $ on $\mathcal{D}$ to approximate $ f(x) $. [\ref{ssec:training}]
    \State Initialize set of quasi-randomly distributed input parameters $ x_{\text{init}} $ across initial dataset. [\ref{ssec:opt}]
	\State Using stochastic gradient descent to find the set $ x_{\text{local}}$ that locally minimizes 	$g(\hat{f}(x), y_\text{target})$. [\ref{ssec:objfun}]
	\State Select the input $ x_{\text{best}}$ that produces minimum loss according to objective function. [\ref{ssec:opt}]
	\State Intelligently query $ f(x) $ with the selected $ x_{\text{best}}$.  [\ref{ssec:query}]
	\If{stopping criteria $= True$}
        \State Take $ x_{\text{best}}$ and $ y_{\text{best}} = f(x_{\text{best}})$.
        \Else
        \State Add $ x_{\text{best}}$ and $ y_{\text{best}} $ to the training dataset $\mathcal{D}$.
        \State Iterate from Step 2. [\ref{ssec:iter}]
    \EndIf  \label{alg:algo1}
\end{algorithmic}
\end{algorithm}
\vspace{-0.15in}

\subsection{Reference Trajectory Planning Algorithm and Simulation} \label{ssec:alip}

Designing the approach and landing trajectory for a hypersonic vehicle is very challenging due to the vehicle's lack of thrusters and large complex physical uncertainties. The process is so computationally expensive that the nominal trajectory is calculated offline and then stored in an onboard computer for the guidance and control system to follow using speed brakes and unboosted maneuvers.

Our simulation uses the Autolanding I-load Program (ALIP) \cite{barton} to calculate the A/L reference trajectory for the HV. This program has legacy from Orbital Science's X-34 and NASA's Shuttle Program \cite{grubler}. ALIP relies on initializing the parameters defining the geometric segments of the trajectory with an accurate prediction and then optimizes those parameters to achieve the desired dynamic pressure on the vehicle at touchdown. Various physical constraints are applied to reduce the trajectory design problem to a two-point boundary value problem \cite{barton}. 

There are 13 input parameters to ALIP that define the initial trajectory that all must be optimized to result in a successful landing of the HV. These parameters include initial dynamic pressure, landing velocity, and flight path angle of constant-glide segments, among others. Fig.~\ref{fig:ALIP} shows an example reference trajectory broken into its geometric segments.

   \begin{figure}[t]
      \centering
      \includegraphics[width=0.99\columnwidth]{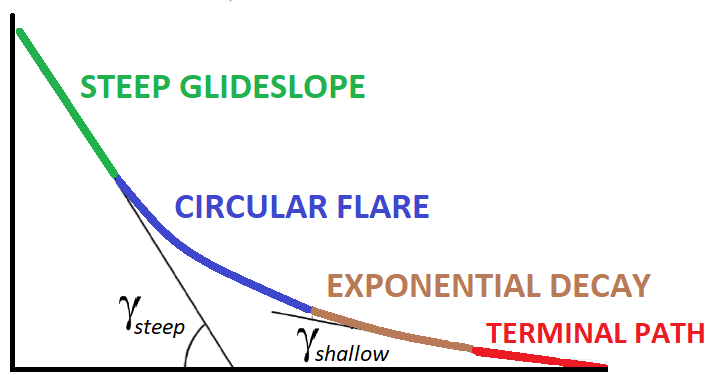}
      \vspace{-0.10in}
      \caption{Defined geometric segments of reference trajectory for approach and landing of hypersonic vehicle.}
      \label{fig:ALIP}
      \vspace{-0.20in}
      \end{figure}

    \begin{figure*}[!b]
    \centering
    \vspace{-0.1in}
    \includegraphics[width=0.8\textwidth]{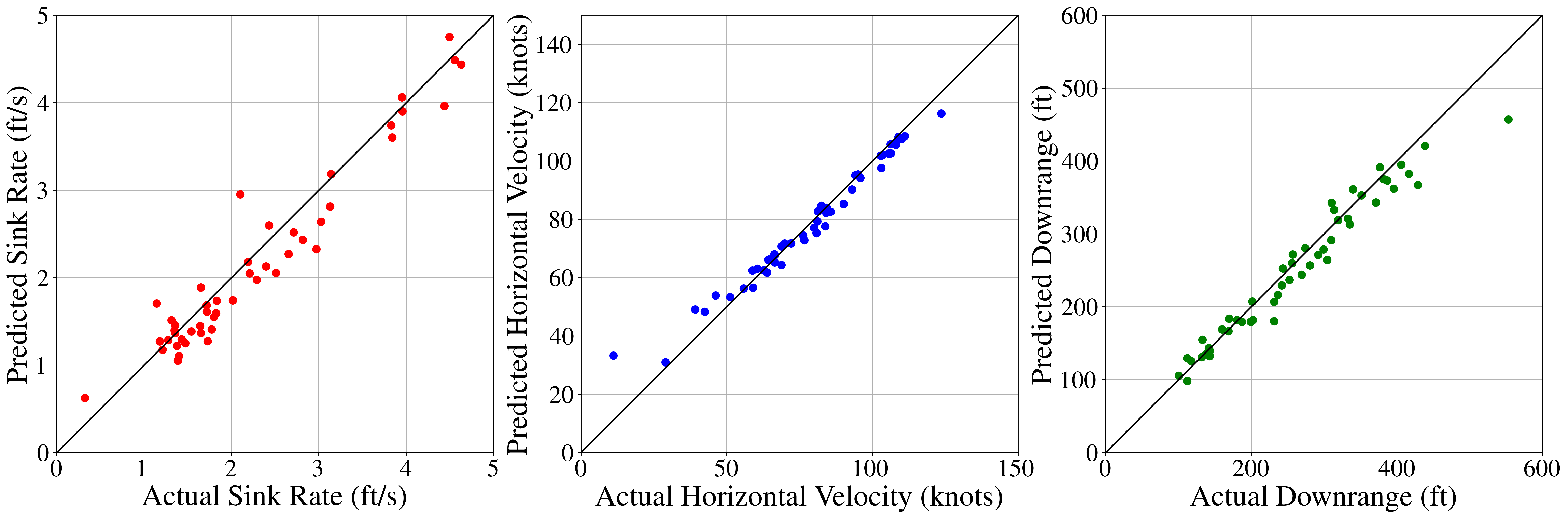}
    \setcounter{figure}{2}
    \vspace{-0.15in}
    \caption{Predicted outputs from surrogate neural network versus actual outputs from simulation. Ideally, these outputs would be perfectly correlated but the neural network is able to mostly accurately predict all three outputs with a few outliers. Sink rate is the most challenging output to predict.}
     \label{fig:output3}
	\end{figure*} 
 
Once a reference trajectory has been defined using ALIP, that trajectory is tested through a digital twin simulation. This simulation outputs certain metrics which quantify the success of the HV touchdown. Specifically, there is an ideal rate at which the vehicle loses altitude (sink rate), distance from the start of the runway where landing happens (downrange position), and horizontal velocity of the HV at touchdown. The digital twin simulation is executed using physical parameters from the hypersonic vehicle. 

All inputs parameters to ALIP and all output parameters from the landing simulation have been standardized and selectively scaled to protect any proprietary information.

\subsection{Initial Data Collection} \label{ssec:initdata}

Each input to the simulation is bounded by physical constraints inside the trajectory planning simulation. These bounds can be adjusted based on historical knowledge of the simulation to reduce the size of the input parameter space. The 13-dimensional input parameters are bounded around a known set of nominal inputs that produce an acceptable landing of the hypersonic vehicle. To create sufficient coverage of the high-dimensional input parameter space, the initial input values were sampled using a Sobol sequence \cite{burhenne}. Sobol sequences are low-discrepancy, quasi-random sequences that provide a more uniform distribution than compared to other quasi-random sampling methods such as Latin hypercube ~\cite{donovan}. 

Enough data samples from the simulation must be collected to sufficiently train the neural network. Too few samples will lead to the surrogate model being a poor representation of the actual simulation and too many samples lead to an unnecessarily high computation cost and time. 

Although the input parameter space is large, the number of samples can be reduced because the surrogate model only needs to understand the gradients of the predicted outputs based on the inputs. An initial sample size of 400 points was chosen to collect input-output mapping for training the surrogate neural network. Increasing the size of the training data past 400 points leads to diminishing improvements to the surrogate's model accuracy, as shown in Fig.~\ref{fig:trainingsize}.

      \begin{figure}[t]
      \centering
      \includegraphics[width=0.99\columnwidth]{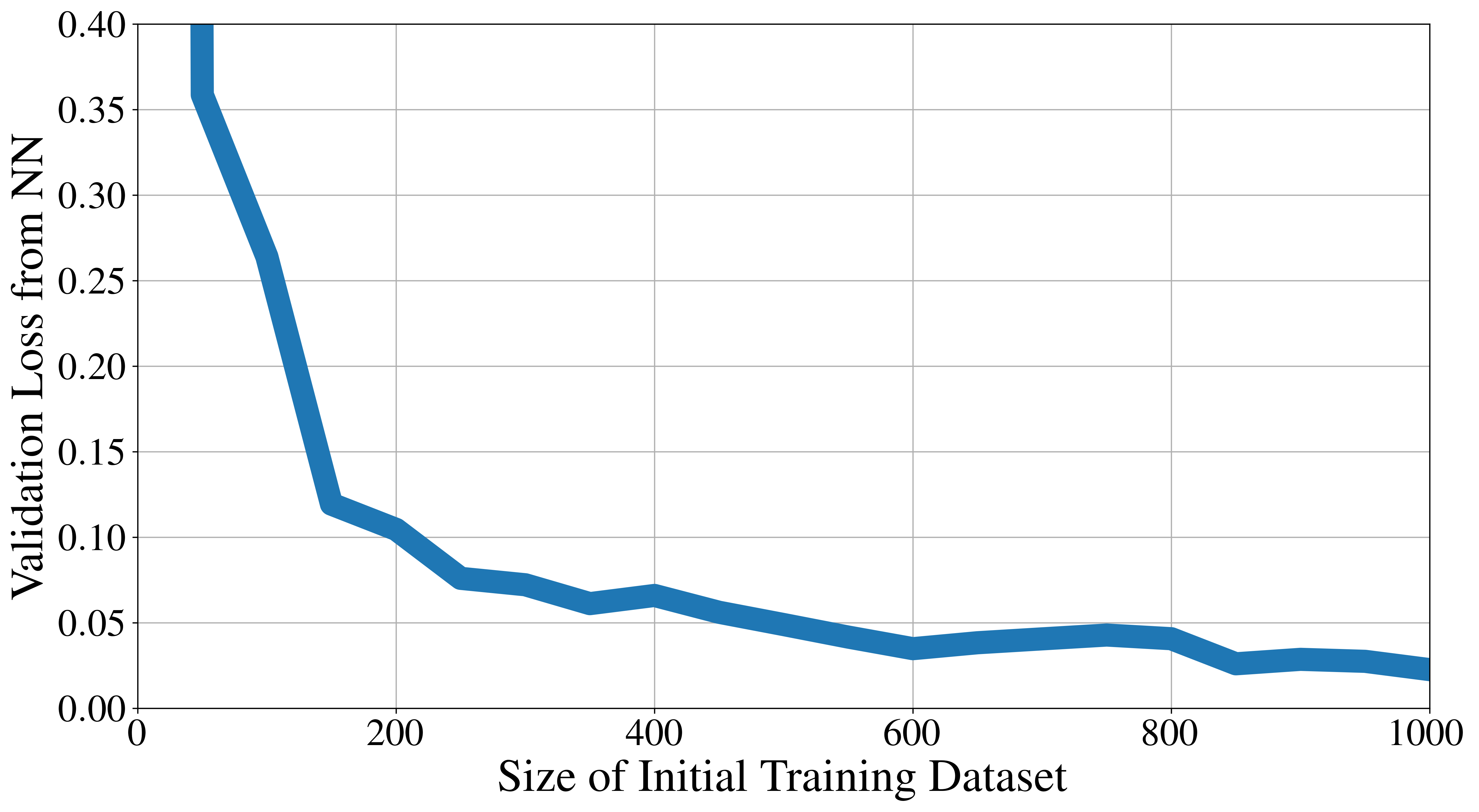}
      \vspace{-0.3in}
      \setcounter{figure}{1}
      \caption{Normalized validation losses from surrogate neural network depending on size of training data set. Diminishing returns on minimizing validation loss after 400 data points.}
      \label{fig:trainingsize}
      \vspace{-0.2in}
   \end{figure}

\subsection{Surrogate Model Training} \label{ssec:training}

Once the initial sampling is done, the data must be standardized for training a neural network to be the best possible surrogate model to accurately represent the simulation. The neural network is trained to minimize the mean absolute error (MAE) between the outputs predicted by the network and the outputs from the training data using the Adam optimizer in PyTorch. The weights of the network are tuned to reduce the MAE using a variant of stochastic gradient descent with backpropagated gradients. The model can be further improved by tuning the hyperparameters defining the neural network's numbers of layers, number of nodes per layer, batch size, and learning rate. The hyperparameters were optimized using the Optuna software library \cite{akiba} integrated with PyTorch. 

Ideally, the neural network would be able to perfectly predict the three outputs of the simulation based on the same set of inputs. Due to the highly nonlinear nature of the simulation and the limited size of the training data, this is not possible. A comparison of the predicted versus actual outputs is shown in Fig.~\ref{fig:output3} for all three outputs of the simulation. 

Fig.~\ref{fig:trainedoutputs} shows the test loss for each output as it converges over epochs. The scale of the test loss varies based on the outputs feasible range in their physical units. For example, the possible range for sink rate only varies between 0 and 5 ft/s, so the surrogate model predicting to an accuracy under 0.2 ft/s is quite good. Similarly, when the feasible range of physical values for downrange position is 0 to 500 feet, predicting within 10 feet is good.

\begin{figure}[t]
      \centering
      \includegraphics[width=0.99\columnwidth]{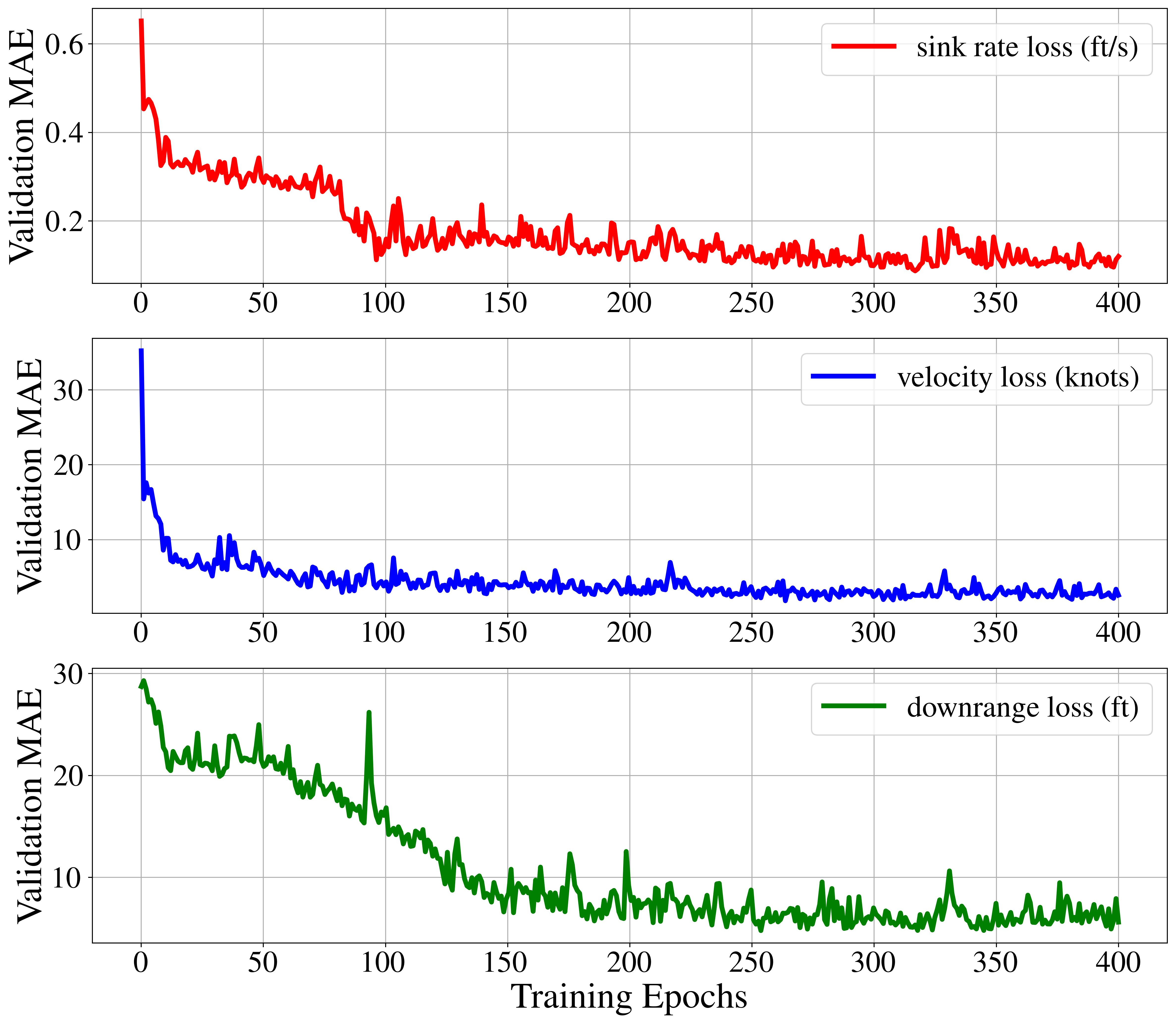}
      \vspace{-0.15in}
      \setcounter{figure}{3}
      \vspace{-0.15in}
      \caption{Validation losses for each simulation output over training epoch. All three outputs converge to low loss relative to their scale.}
      \label{fig:trainedoutputs}
      \vspace{-0.2in}
   \end{figure}

\subsection{Custom Objective Function} \label{ssec:objfun}

The objective of this constrained optimization problem is to minimize the distance between desired and predicted outputs. For this application, the objective function to be minimized evaluates the mean absolute error between the desired outputs and the predicted outputs from the surrogate model. MAE calculates the average of the absolute distance between all the predicted and desired outputs. MAE is more robust to outliers and ensure stable gradients with the different terms in the objective function. This makes the objective function
\begin{equation} g(\hat{y})=\left| y_{\text {target }} - \hat{y} \right|. \end{equation}
However, the parameters being optimized are not the three predicted outputs from the surrogate model $\hat{y}$, but the 13 inputs to the surrogate model $x$, making the objective function
\begin{equation} g(x)= \left| y_{\text {target }} - \hat{f}(x)\right|. \end{equation}
Furthermore, each of the desired outputs can be weighted based on domain knowledge of the simulation. In this application, sink rate is more important to a successful landing than either horizontal velocity or downrange position. Therefore, the MAE of predicted versus actual sink rate is scaled by a factor of two. And the desired downrange position is too high to be realistically achievable by any combination of input values, so its MAE is reduced by a factor of 10 to ensure it does not overpower the other two terms.

To constrain the inputs to the physically feasible parameter space, the inequality constraints are added directly to the objective function as a soft penalty. When either the upper or lower bounds are violated, these penalties are triggered as
\begin{equation}
  \begin{gathered}
	g(x)=\left| y_{\text {target}}- \hat{f}( x )\right|+\alpha \max \left(0, x_{\min }- x \right) \\
	+ \alpha \max \left(0, x - x_{\max }\right) \label{eqn:objective}
  \end{gathered}
\end{equation}
where $\alpha=1$, but can be tuned depending on the scale of the other components of the objective function. These penalties act to steer the inputs being optimized to stay in the constrained input parameter space through gradient descent. The gradients of the objective function with respect to the inputs are used to optimize the inputs to minimize the objective. The gradient of the objective function with respect to the inputs when the input bounds are violated is 
\begin{equation} \frac{\partial g}{\partial x}=\hat{f}^{\prime}(x)\left(\frac{\hat{f}(x)-y_{\text {target }}}{\left| y_{\text {target }}-\hat{f}(x) \right| } \right) \pm \alpha \end{equation}
As defined in \ref{ssec:alip}, the outputs are sink rate, downrange position, and horizontal velocity of the HV at touchdown. The target values of the outputs is determined by domain knowledge and vehicle requirements. It is important to note that the number of desired outputs as well as their values could be changed, and this novel algorithm would still work. For example, if the requirement for sink rate was lowered to 1.5 ft/s, the objective function can  be adjusted to produce a new optimal reference trajectory. 

\subsection{Sensitivity Analysis and Loss Landscape Visualization}

\begin{figure}[t]
      \centering
      \includegraphics[width=0.99\columnwidth]{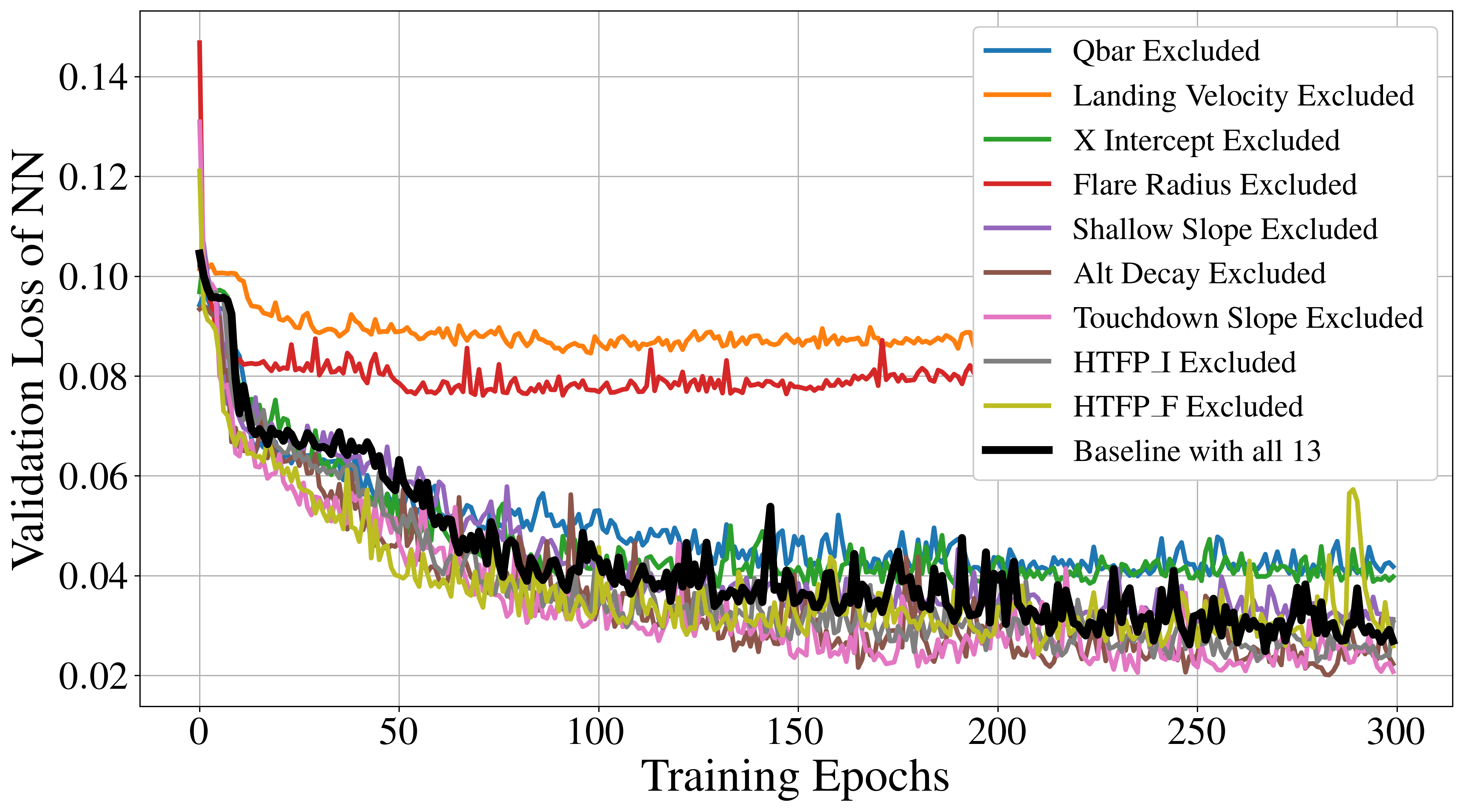}
      \vspace{-0.3in}
      \caption{Sensitivity Analysis Results. Landing velocity and flare radius inputs have the most significant impact on the accuracy of the surrogate model.}
      \label{fig:sensitivity}
      \vspace{-0.2in}
   \end{figure}

It is impossible to easily visualize the gradients of the objective function with respect to all 13 inputs. While it would be possible to do a dimensional reduction to two dimensions, because these inputs have real meanings the two most important inputs can simply be identified. A sensitivity analysis was conducted to identify the two most important inputs, meaning the two inputs that most affect the output of the simulation. This was done by excluding one input at a time from the NN training and evaluating which input increased the loss the most, therefore most effecting the neural network's ability to accurately predict the outputs. In Fig.~\ref{fig:sensitivity}, the two most important inputs for predicting the outputs are shown to be the landing velocity of the HV and the circular flare radius.

\begin{figure}[t]
      \centering
      \includegraphics[width=0.99\columnwidth]{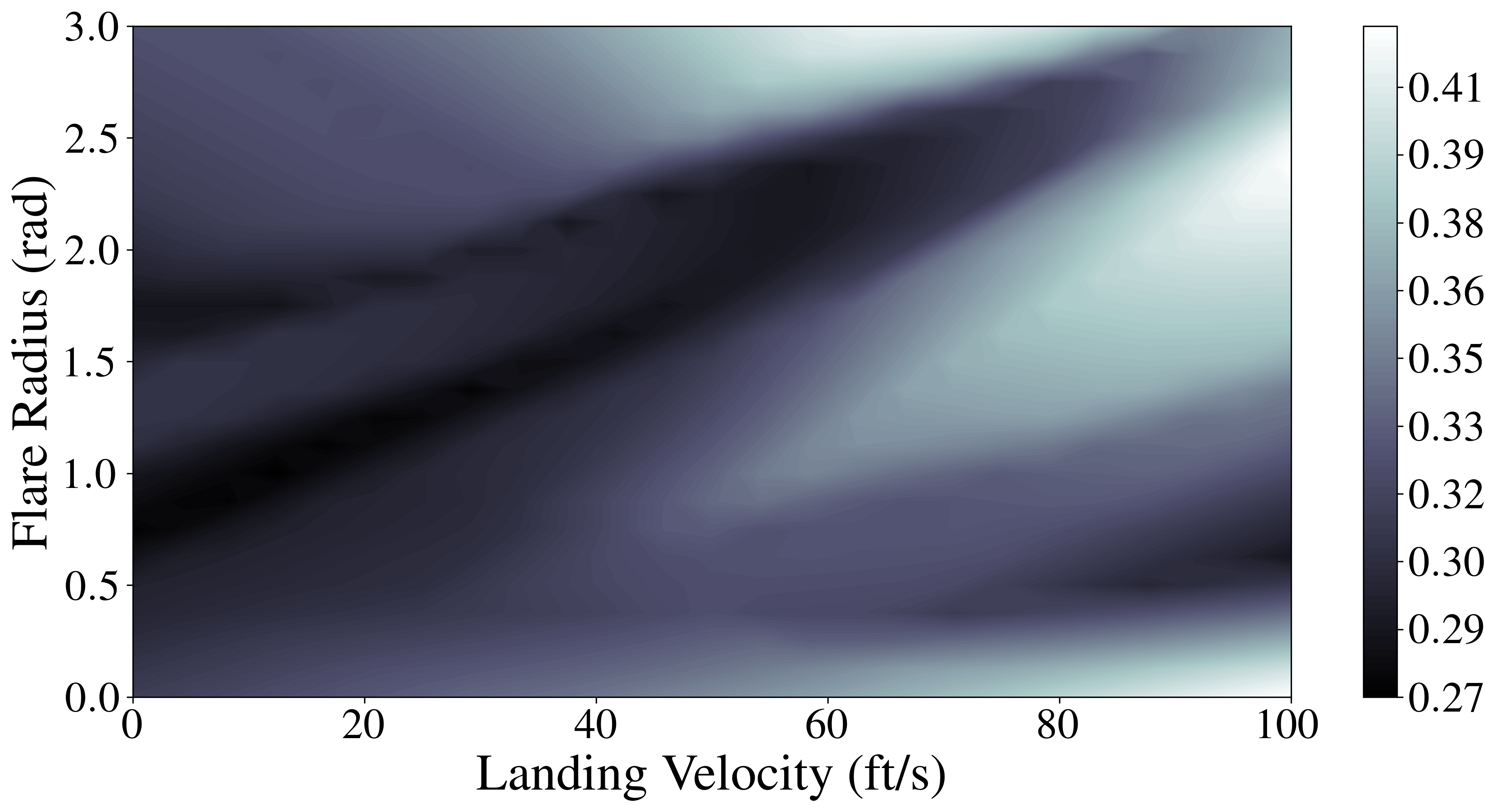}
      \vspace{-0.3in}
      \caption{Loss landscape as predicted by surrogate neural network.}
      \label{fig:predloss}
      \vspace{0.1in}
      \centering
      \includegraphics[width=0.99\columnwidth]{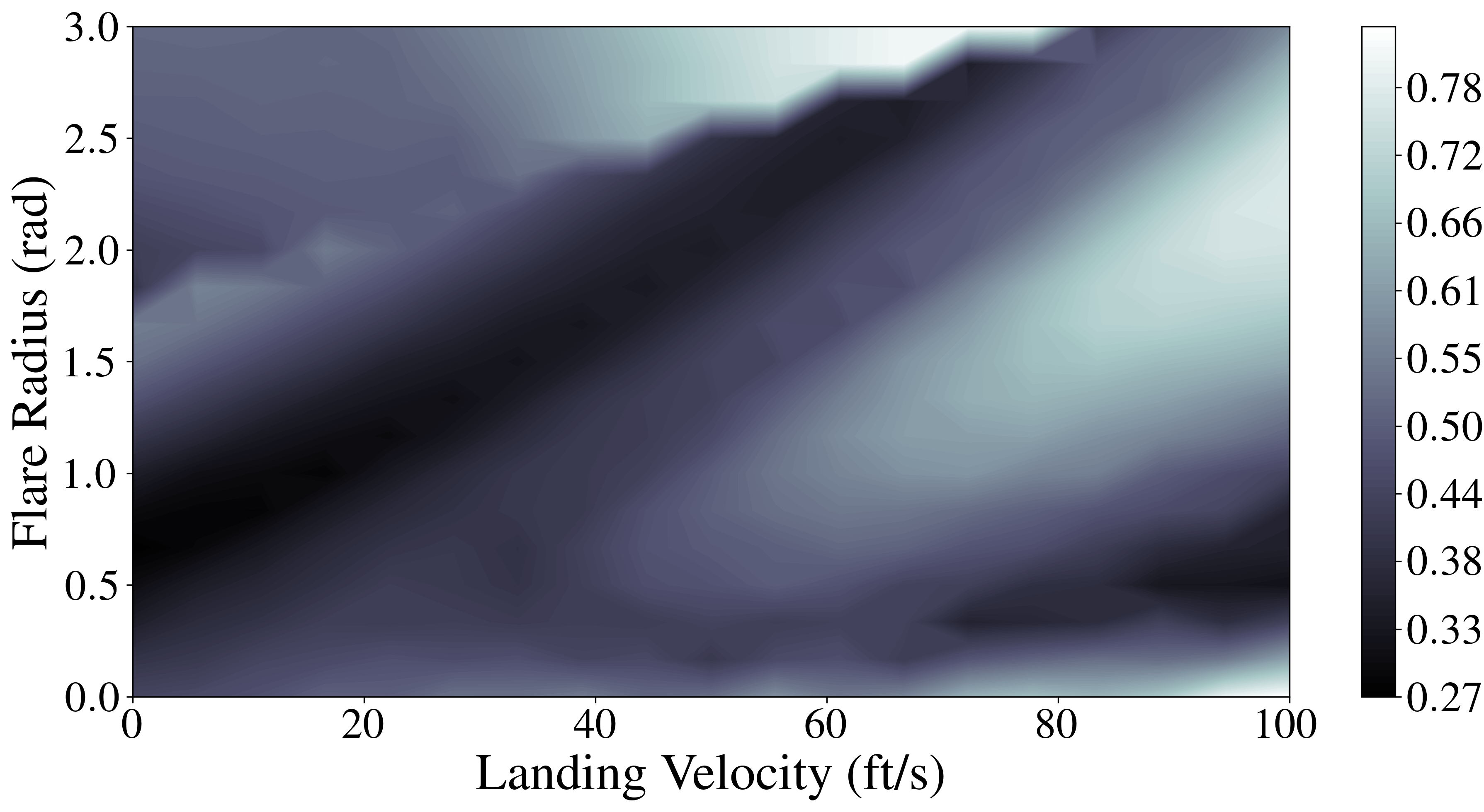}
      \vspace{-0.3in}
      \caption{Loss landscape from real outputs of the simulation. The loss scale and gradients match those of the trained surrogate model.}
      \label{fig:realloss}
      \vspace{-0.25in}
   \end{figure}

With the two most important inputs identified, how varying the inputs affect the objective function can be visualized. While not necessary for optimization, these visualizations are to better understand and validate this methodology. Linearly spaced samples across the input range for landing velocity and flare radius were selected while all other inputs were frozen to their nominal values. Each combination of landing velocity and flare radius was fed through the trained surrogate neural network, and the corresponding output was evaluated by the objective function to calculate loss. This loss landscape in Fig.~\ref{fig:predloss} shows how the loss from the objective function varies across the two most important inputs, as predicted by the surrogate model.
   
It is evident that the surrogate neural network has identified regions in the input parameter space as more promising for minimizing the objective function than others. Specifically, a trough in the middle where both variables increase proportionally and a section of low flare radius and high landing velocity. To validate the surrogate model, the true loss landscape from actually querying the simulation is shown in Fig.~\ref{fig:realloss}.

\subsection{Input Optimization and Intelligent Querying} \label{ssec:opt}

With the initial surrogate model trained, the 13 inputs can now be optimized to minimize the objective function using gradient-based descent. Similar to how the gradients of the neural network are used to optimize the weights during the initial training, the gradients are used to optimize the inputs to the simulation to minimize the objective function (\ref{eqn:objective}). 

One of the issues with gradient-based optimization is the convergence to local minima or maxima. To mitigate this concern, many input vectors are optimized according to the objective function. These input vectors are quasi-randomly distributed across the input parameter space using a Sobol sequence.  Even if some of the inputs become stuck at local minima or plateaus, at least one of the vectors will find the global minimum. 

Through trial-and-error, the best optimizer was found to be stochastic gradient descent with added momentum. Because the many inputs are optimized in batches, this process is extremely rapid compared to the time it takes to query the simulation, seconds as opposed to minutes. This means it is computationally cheap to initialize hundreds of inputs and optimize over many steps to ensure convergence. In the methodology proposed in this paper, 100 input vectors are optimized over 500 steps. 

\begin{figure}[t]
      \centering
      \includegraphics[width=0.99\columnwidth]{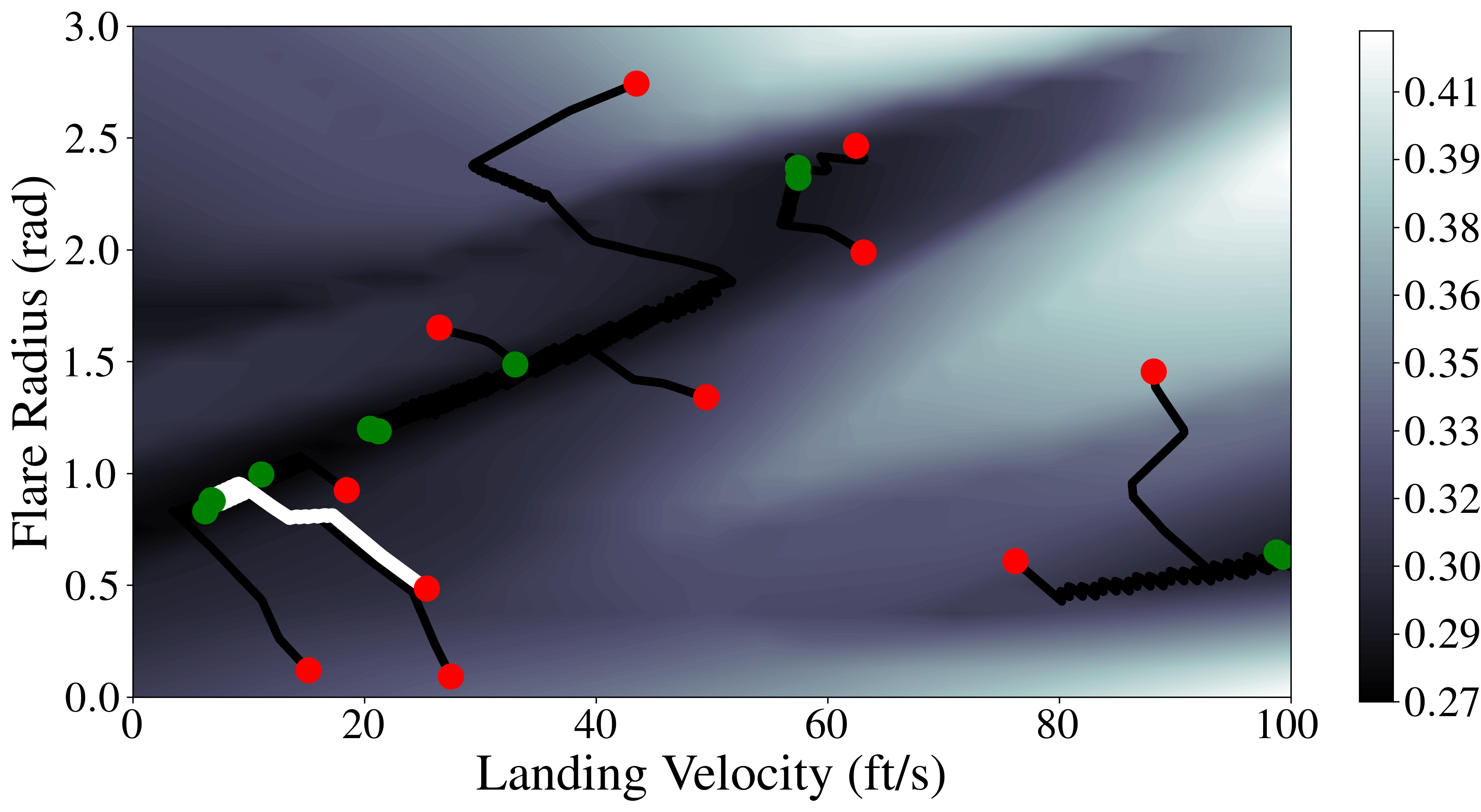}
      \vspace{-0.3in}
      \caption{Optimization of inputs across loss landscape. The input vectors are optimized from places with higher loss to areas that, according to the surrogate model, will reduce the loss from objective function.}
      \label{fig:optimize}
      \vspace{-0.2in}
   \end{figure}

This process is visualized in Fig.~\ref{fig:optimize} using the same loss landscape from the previous section. Again, the inputs are limited to the two most important parameters, landing velocity and flare radius, so the optimization corresponds to the loss landscape background. Only 10 of the 100 input vectors are shown, their initial and final points shown in red and green, respectively. The input vector that finds the global minimum has its path in white instead of black.

The input vectors clearly travel down the gradient from places with higher estimated loss to places with lower estimated loss over optimization. It is also shown that the input vectors are successfully constrained within the bounds.

\subsection{Intelligent Querying} \label{ssec:query}

The input vector that found the global minimum of the loss landscape is now sent to the actual simulation. While there are many input vectors that show promising results for minimizing the objective function according to the surrogate model, only one is sent to the trajectory planning simulation. The input selected to query the simulation is intelligently selected instead of being randomly selected or intuitively chosen by an engineer. This reduces the overall computational cost of optimizing the simulation and finding an exceptional reference trajectory for the HV.

\subsection{Iteration} \label{ssec:iter}

Until the stopping criteria is met, this algorithm will iteratively train a new neural network, optimize the objective function, and intelligently query the simulation. The stopping criteria for this application checks for three conditions every iteration:

\begin{enumerate}

\item If the goal criteria for the outputs have been met
\item If the algorithm converged on a solution outside of the goal criteria
\item If the maximum number of iterations been met 

\end{enumerate}

Each iteration of the algorithm allows the surrogate model to improve its representation of the actual system and increase its fidelity in areas that appear promising for minimizing loss. Without initializing and training a new neural network each time, the algorithm would always search around the same global minimum. Instead, the random weight initialization for each iteration's surrogate model means each neural network will be slightly different. 

Furthermore, having only one surrogate model does not capture any `model uncertainty'. Creating a new model every time is more similar to ensemble learning and is shown to be more robust to outliers \cite{laksh}. This encourages more exploration during optimization instead of continually exploiting the same promising regions. 

\section{NUMERICAL RESULTS}

\subsection{Optimal Reference Trajectory} \label{ssec:traj1}

This section demonstrates the success of this algorithm when applied to optimizing the reference trajectory for the hypersonic vehicle. This algorithm found a new reference trajectory for approach and landing that outperformed the existing nominal solution. The new reference trajectory produced a 74\% decrease in loss from the objective function when compared to the nominal solution as shown in Table \ref{tab:perfcomp}.

\begin{table}[b]
\vspace{-0.15in}
\caption{FIRST REFERENCE TRAJECTORY COMPARISON}
\vspace{-0.15in}
\label{tab:perfcomp}
\begin{center}
\begin{tabular}{|l|c|c|c|c|}
\hline \multirow{2}{*}{} & \multicolumn{4}{c|}{ Performance Results of Trajectories } \\
\cline { 2 - 5 } & 
$\begin{gathered}\text{ Sink Rate } \\[-1mm] (\text{ft/s) }\end{gathered}$ & 
$\begin{gathered}\text{ Horizontal } \\[-1mm] \text{ Velocity }(\text {knots) }\end{gathered}$ & $\begin{gathered}\text{ Downrange } \\[-1mm] \text{ Position }(\text {ft})\end{gathered}$ & \multicolumn{1}{c|}{ Loss } \\\hline 
Target & $2.00$ & $54.0$ & $400.0$ & $0.00$ \\
Nominal & $1.48$ & $69.9$ & $193.2$ & $0.24$ \\
New & $2.00$ & $53.8$ & $378.4$ & $0.11$ \\
\hline
\end{tabular}
\end{center}
\end{table}

The algorithm is very successful at finding values for the 13 simulation inputs that result in values for sink rate and horizontal velocity almost exactly at their target values. However, downrange position is harder to achieve due to internal simulation constraints. The algorithm is essentially trying to maximize downrange position without compromising the results for sink rate and horizontal velocity.  

\begin{figure}[t]
      \centering
      \includegraphics[width=0.99\columnwidth]{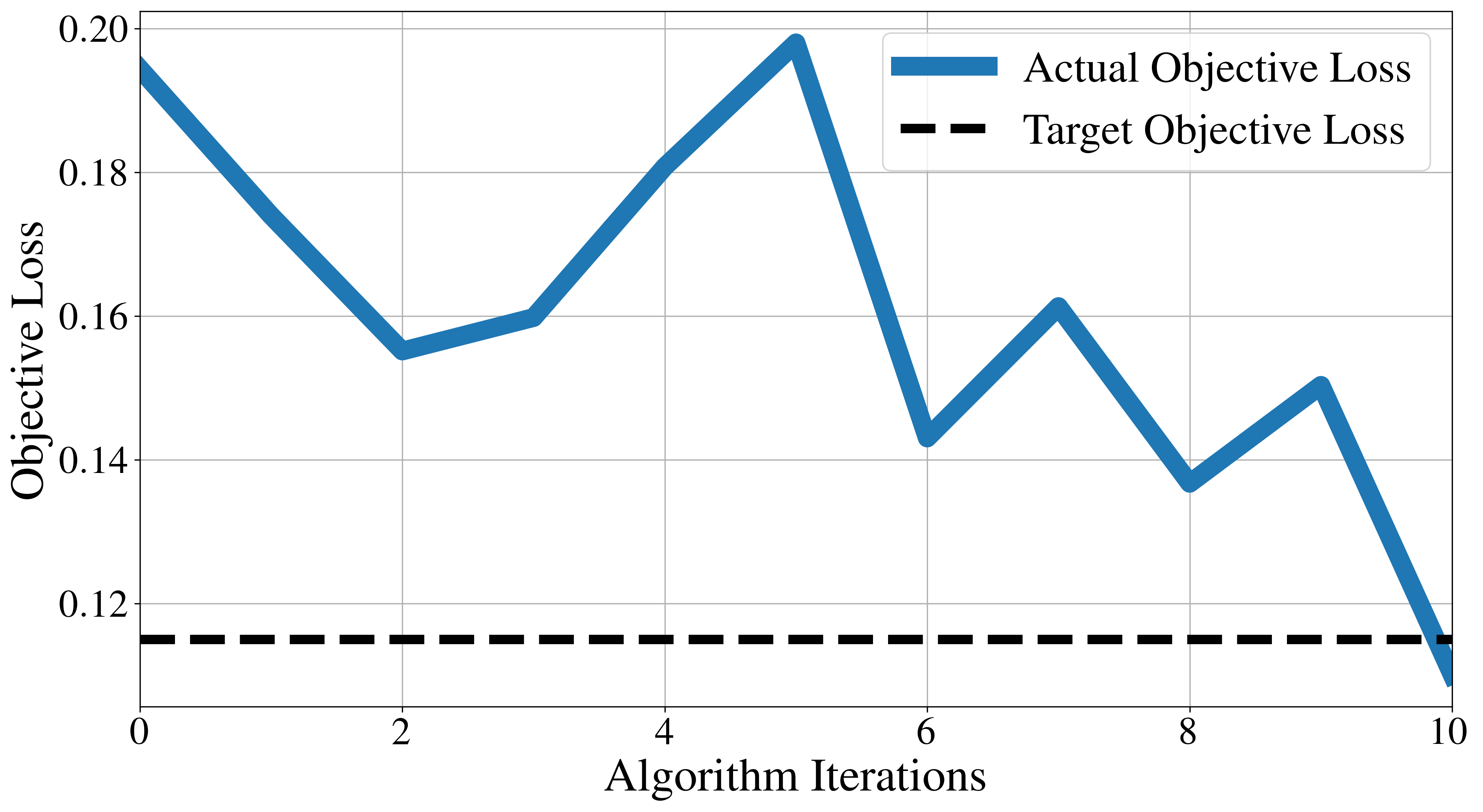}
      \vspace{-0.3in}
      \caption{Objective loss over algorithm iterations until successful simulation outputs are achieved.}
      \label{fig:bestloss}
      \vspace{-0.2in}
   \end{figure}

This new optimal reference trajectory was obtained by training the surrogate NN on an initial dataset of 400 points quasi-randomly selected across the input parameter space. The set of input parameter values that produced this optimal reference trajectory were found after only 10 queries to the actual simulation. The objective loss of the true outputs from the simulation is shown over those 10 iterations until the end condition of reaching desired values is satisfied in Fig.~\ref{fig:bestloss}.

\subsection{Monte Carlo Comparison}

Although this methodology has proved its ability to optimize the trajectory planning simulation, the true motivation of this algorithm is in its timesaving. It is able to intelligently and efficiently search to the input parameter space to find an optimal solution faster than tuning by hand or through a Monte Carlo random search approach.

   \begin{figure}[t]
      \centering
      \includegraphics[width=0.99\columnwidth]{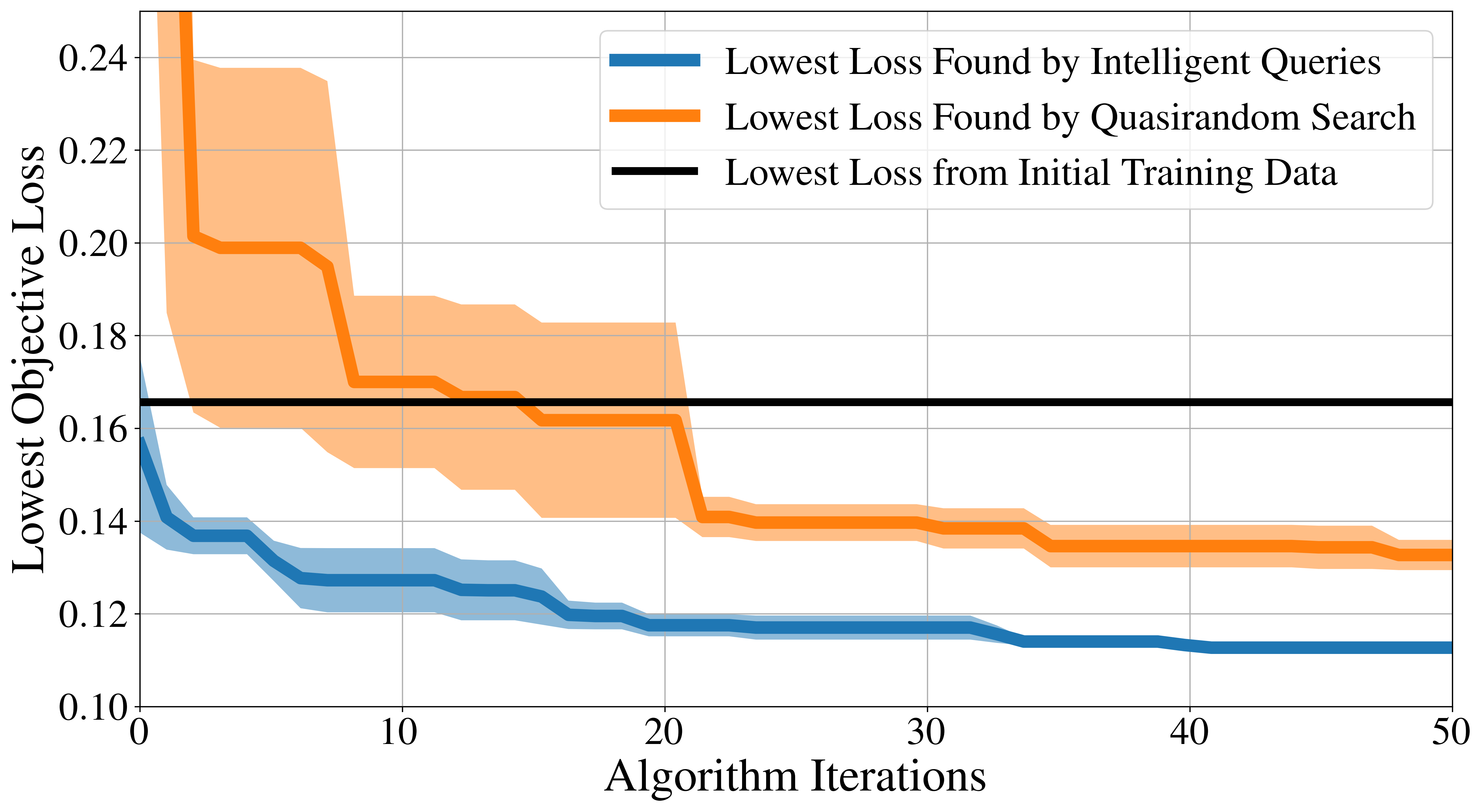}
      \vspace{-0.3in}
      \caption{Average lowest loss from intelligent queries compared to quasi-random queries over 50 algorithm iterations with 0.5$ \sigma $ shaded. The intelligent queries result in lower losses quicker than quasi-random queries.}
      \label{fig:lowestloss}
	\vspace{0.1in}
      \includegraphics[scale=0.22]{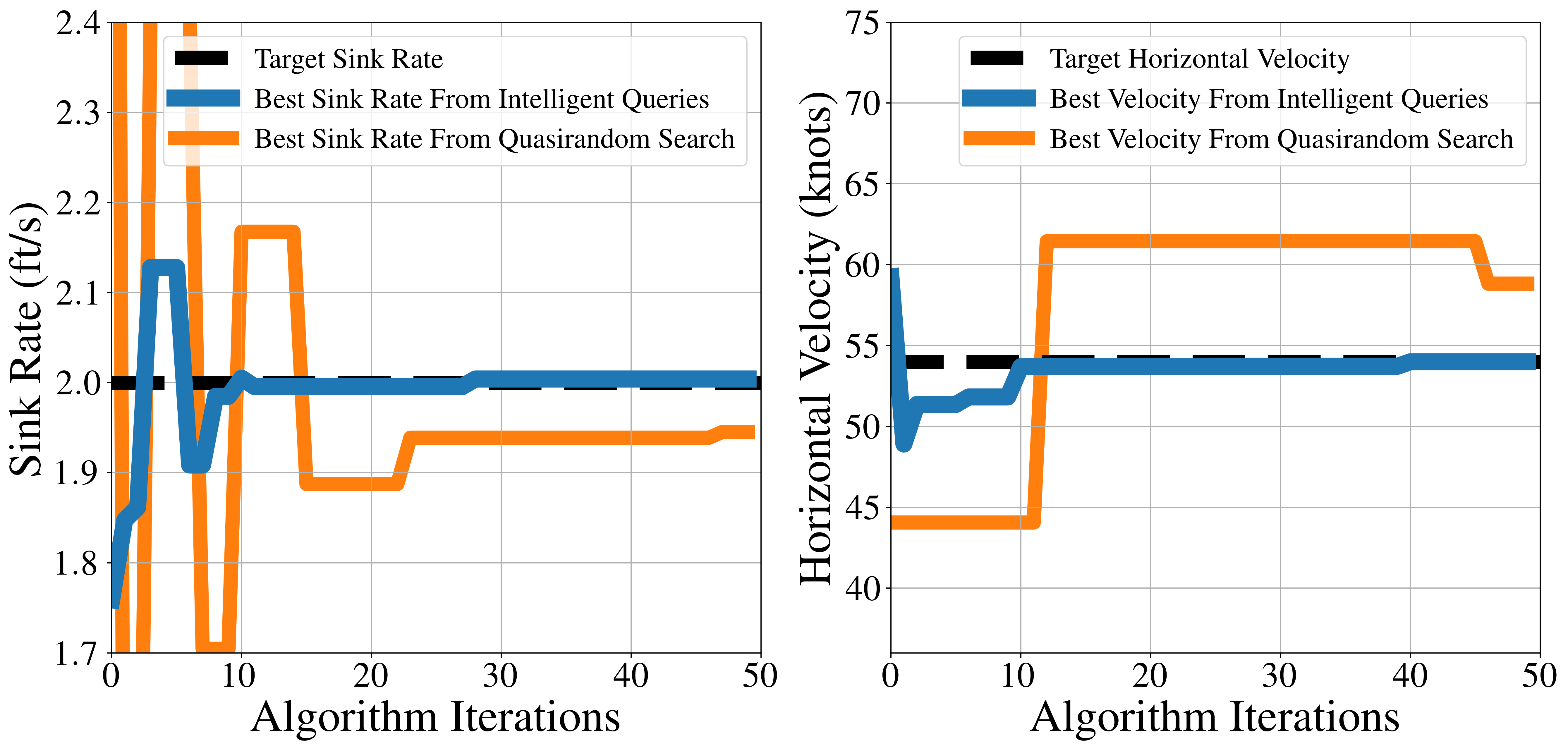}
      \vspace{-0.3in}
      \caption{Best outputs found by intelligent queries compared to Monte Carlo quasi-random search. The intelligent queries converge much faster to the ideal outputs than the quasi-random search.}
      \label{fig:closestouts}
      \vspace{-0.2in}
   \end{figure} 
   
This can further be proved by tracking the lowest objective loss found so far by the quasi-random sampling and this intelligent algorithm. Over 50 iterations, this approach consistently finds better solutions faster with less variation than quasi-random sampling. Fig.~\ref{fig:lowestloss} shows the loss from the best solution found so far by this algorithm as compared to the quasi-random search, averaged over five trials.
   
Evidently, this algorithm finds inputs that minimize the objective function much faster than a quasi-random search can. The quasi-random search takes 30 simulation queries to reduce the objective loss from 0.16 to 0.13, while the intelligent search reduces the same loss in the first five queries, making this method six times faster. Furthermore, it consistently reaches a lower objective loss, meaning it finds a better reference trajectory. Similarly, Fig.~\ref{fig:closestouts} shows how the algorithm will converge to the desired output values faster than the quasi-random search.

\begin{table}[b]
\vspace{-0.15in}
\caption{SECOND REFERENCE TRAJECTORY COMPARISON}
\vspace{-0.15in}
\label{tab:perfcomp2}
\begin{center}
\begin{tabular}{|l|c|c|c|c|}
\hline \multirow{2}{*}{} & \multicolumn{4}{c|}{ Performance Results of Trajectories } \\
\cline { 2 - 5 } & $\begin{gathered}\text { Sink Rate } \\[-1mm]
(\text {ft/s) }\end{gathered}$ & $\begin{gathered}\text { Horizontal } \\[-1mm]
\text { Velocity }(\text {knots) }\end{gathered}$ & $\begin{gathered}\text { Downrange } \\[-1mm]
\text { Position }(\text {ft})\end{gathered}$ & \multicolumn{1}{c|}{ Loss } \\\hline 
Target & $2.00$ & $54.0$ & $400.0$ & $0.00$ \\
Previous & $1.21$ & $52.5$ & $365.8$ & $0.51$ \\
New & $1.76$ & $53.4$ & $354.0$ & $0.17$ \\
\hline
\end{tabular}
\end{center}
\end{table}

\subsection{Second Hypersonic Vehicle Simulation}

All the results so far show the optimization of an A/L reference trajectory based on a simulation of a proprietary hypersonic vehicle. If the underlying physical parameters of the vehicle change, a new reference trajectory will need to be calculated to meet the new set of requirements. For example, if additional wind tunnel testing results in an update to the vehicle aerodynamic model. This methodology can quickly be applied to calculate a new optimal reference trajectory, negating the need for an engineer to spend copious time tuning the input variables. 

To test this, the existing HV's aerodynamic properties were perturbed randomly within a range of 3$\sigma$. Specifically, the drag was increased by 11\%, the lift was reduced by 1\% and the pitch moment coefficient was increased by 13\%. Now, the set of inputs found in Subsection \ref{ssec:traj1} result in a very low sink rate and therefore a poor landing of the HV. When applied, this methodology finds a better performing reference trajectory that decreases loss by 100\% as shown in Table \ref{tab:perfcomp2}.

\section{CONCLUSIONS}

This paper presents a new simulation-based optimization algorithm that uses surrogate neural networks and their ability to produce gradients to drastically reduce computation time. When applied to a highly nonlinear A/L reference trajectory planning simulation for hypersonic vehicles, this algorithm rapidly optimized the 13 input parameters to produce the best possible result for landing the hypersonic vehicle. The novelty of this algorithm is that it uses neural network surrogate models to intelligently select queries to the simulation. By doing so, the total number of simulation runs can be reduced while simultaneously finding the optimal reference trajectory. 

This generalized methodology has been shown to work for different HV simulations. However, future work includes testing this methodology in an entirely different field. This methodology should be able to accommodate any black-box simulation with any number of inputs or outputs. Furthermore, total number of queries could be reduced by eliminating the initial quasi-random training dataset so that all queries are intelligently selected by the surrogate model. It would also be interesting to explicitly incorporate the epistemic uncertainty of the neural network to better balance and understand exploration versus exploitation.

Compared to quasi-random search, this algorithm works six times faster to find the optimal output of the simulation by intelligently querying areas that minimize loss as predicted by the surrogate model. Furthermore, no current optimal trajectory planning algorithms make use of surrogate models, their efficiency, and their ability to produce gradients. This methodology will enable much more rapid calculation of optimal approach and landing reference trajectories for hypersonic vehicles.

\balance

\addtolength{\textheight}{-1cm}   


\end{document}